\documentclass[11pt]{article}
\usepackage{bcc24}

\bcctitle{Some open problems on permutation patterns}
\bccshorttitle{Some open problems on permutation patterns}

\bccname{Einar Steingr\'{\i}msson\footnote{Supported by grant no.\
    090038013 from the Icelandic Research Fund.}}
\bccshortname{E. Steingr\'{\i}msson}

\bccaddressa{Department of Computer and Information Sciences\\University
  of Strathclyde\\Glasgow G1 1XH, UK.}
\bccemaila{einar@alum.mit.edu}
\usepackage[utf8]{inputenc}
\usepackage[T1]{fontenc}

\usepackage{tikz}
\usetikzlibrary{matrix,arrows,patterns}

\DeclareMathOperator{\asc}{\mathrm{asc}}
\DeclareMathOperator{\sgn}{\mathrm{sgn}}

\newcommand{\clp}{\mathcal{P}}
\newcommand{\cli}{\mathcal{I}}

\newcommand{\clc}{\mathcal{C}}
\newcommand{\cln}{\clc_n}

\newcommand{\cls}{\mathcal{S}}

\newcommand\avpn{\ensuremath{\mathrm{Av}_p(n)}}

\newcommand\avsn{\ensuremath{\mathrm{Av}_\sigma(n)}}
\newcommand\avtn{\ensuremath{\mathrm{Av}_\tau(n)}}

\newcommand\ist{\ensuremath{[\sigma,\tau]}}
\newcommand\mst{\ensuremath{\mu(\sigma,\tau)}}
\newcommand{\one}{\mathbf{1}}  %the permutation of length 1

\def\ch#1,#2,{{\binom{#1}{#2}}}

\def\bskip#1,{\vspace*{#1\baselineskip}}
\def\emm#1,{{\em#1}}
\def\cd{\cdot}

\newcommand\op{\oplus} 
\newcommand\om{\ominus}

\newcommand\sw{\mbox{\sc\small{SW}}}

\def\dd{\kern.1ex\mbox{\raise.7ex\hbox{{\rule{.35em}{.15ex}}}}\kern.1ex}

\usepackage{epic, eepic}

\def\xa{0}
\def\xb{50}
\def\xc{85}

\def\ya{0}
\def\yb{30}
\def\yc{60}
\def\yd{90}

\def\yao{5}
\def\ydu{85}

\def\ybo{35}
\def\yco{65}

\def\ybu{25}
\def\ycu{55}
\def\p{\circle*{2}}

\def\cbox#1{\makebox(0,0){#1}}

\usepackage{tikz}
\usetikzlibrary{matrix,arrows,patterns}

\newcommand{\pattern}[4]{
 \raisebox{0.6ex}{
 \begin{tikzpicture}[scale=0.35, baseline=(current bounding box.center), #1]
  \foreach \x/\y in {#4}
    \fill[pattern=north east lines, pattern color=black!45] (\x,\y)
rectangle +(1,1);
  \draw (0.01,0.01) grid (#2+0.99,#2+0.99);
  \foreach \x/\y in {#3}
    \filldraw (\x,\y) circle (5pt);
 \end{tikzpicture}}\;
}

\begin{document}
\makebcctitle

\begin{abstract}
  This is a brief survey of some open problems on permutation
  patterns, with an emphasis on subjects not covered in the recent
  book by Kitaev, \emph{Patterns in Permutations and words}.  I first
  survey recent developments on the enumeration and asymptotics of the
  pattern 1324, the last pattern of length 4 whose asymptotic growth
  is unknown, and related issues such as upper bounds for the number
  of avoiders of any pattern of length $k$ for any given $k$. Other
  subjects treated are the M\"obius function, topological properties
  and other algebraic aspects of the poset of permutations, ordered by
  containment, and also the study of growth rates of permutation
  classes, which are containment closed subsets of this poset.
\end{abstract}

\thispagestyle{empty}

\section{Introduction}\label{sec-intro}

The notion of permutation patterns is implicit in the literature a
long way back, which is no surprise given that permutations are a
natural object in many branches of mathematics, and because patterns
of various sorts are ubiquitous in any study of discrete objects.  In
recent decades the study of permutation patterns has become a
discipline in its own right, with hundreds of published papers.  This
rapid development has not only led to myriad new results, but also,
and more interestingly, spawned several different research directions
in the last few years.  Also, many connections have been discovered
between permutation patterns and other research areas, both inside and
outside of combinatorics, showcasing the fundamental nature of
patterns in permutations and other kinds of words.

Recently, Sergey Kitaev published a comprehensive reference work
entitled \emph{Patterns in Permutations and words} \cite{kitaev-book}.
In the present paper I highlight some aspects of some of the most
recent developments and some areas that have hardly been touched, but
which I think deserve more attention.  This is speculative, of course,
and strongly coloured by my own preferences and by the limits of my
own knowledge in the field.  The topics dealt with here are mostly
left out in \cite{kitaev-book} (due to their very recent appearance)
so there is little overlap here with that book.

In Section \ref{sec-kinds} I briefly describe the different kinds of
patterns that are prominent in the field.  In Section~\ref{sec-1324} I
treat one case of pattern avoidance, that of the pattern 1324, which
is the shortest classical pattern whose avoidance is not known.  Not
even the asymptotics of the number of 1324-avoiders is known, even
though significant effort has been put into this, resulting in ever
better bounds.  In fact, after a hiatus of a few years there have been
several successively improved results in the last year, which also
seem likely to be applicable in a wider context.

Section~\ref{sec-layered} describes conjectures about which patterns
are easiest to avoid, that is, are avoided by more permutations than
other patterns of the same length.  All the evidence points in the
same direction, namely, that for patterns of any given length~$k$ it
is a layered pattern that is avoided by most permutations.  It is
still unclear, however, what form the most avoided layered pattern of
length $k$ has, although there are conjectures for particular families
of values of $k$ that seem reasonable, while others (in an abundance
of conjectures, published and unpublished) have been shown false.

In Section \ref{sec-poset} I discuss the poset (partially ordered set)
$\clp$ consisting of all permutations, ordered by pattern containment.
This poset is the underlying object of all studies of pattern
avoidance and containment.  I mention the results so far on the
M\"obius function of $\clp$, perhaps the most important combinatorial
invariant of a poset, and some topological aspects of the order
complexes of intervals in $\clp$.  Hardly anything is known so far on
the topology of these intervals, but there are indications that large
classes of them have a nice topological structure, whose understanding
might shed light on various pattern problems.

In Section \ref{sec-algebra} I mention some algebraic properties of
the set of mesh patterns, a generalisation encompassing all the kinds
of patterns considered here, and also the ring of functions counting
occurrences of vincular patterns in permutations.  Neither of these
aspects has been studied substantially, but there are reasons to
believe that they might be interesting.

Finally, in Section~\ref{sec-growth} I treat permutation classes,
which are classes of permutations avoiding a set, finite or infinite,
of classical patterns.  The emphasis here is on the growth rates of
these classes, that is, the growth of the number of permutations of
length $n$ in a given class.  Substantial progress has been made here
in recent years, although this only begins to scratch the surface of
what promises to be an interesting study of fundamental properties of
pattern containment and avoidance.  In particular, the study of
permutation classes and their growth rates is intimately related to
the types of generating functions enumerating these classes, and
generating functions are among the most important tools in the study
of permutation patterns, as in all of enumerative combinatorics.

\section{Kinds of patterns}\label{sec-kinds}

We write permutations in one-line notation, as $a_1a_2\ldots a_n$,
where the $a_i$ are precisely the integers in $[n]=\{1,2,\ldots,n\}$.
For terminology not defined here see \cite{kitaev-book}.

An occurrence of a classical pattern $p=p_1p_2\ldots p_k$ in a
permutation $\pi=a_1a_2\ldots a_n$ is a subsequence
$a_{i_1}a_{i_2}\ldots a_{i_k}$ of $\pi$ whose letters appear in the
same relative order of size as those in $p$.  For example,
$1\dd3\dd2$\footnote{In this section I write classical patterns with
  dashes between all pairs of adjacent letters, to distinguish them
  from other vincular and bivincular patterns.  In later sections, I
  will write classical patterns in the classical way, without any
  dashes, to keep the notation less cumbersome.} appears in 31542 as
142, 152, 154 and 354.  Here, 354 is also an occurrence of the
\emph{vincular} pattern $1\dd32$, because the 5 and 4 are adjacent in
31542, which is required by the absence of a dash between 3 and 2 in
$1\dd32$.  Also, 142 is an occurrence of the \emph{bivincular} pattern
$\bar1\dd32$, because 4 and 2 are adjacent (in position) and 1 and 2
are adjacent in value, as required by the bar over the~1 in
$\bar1\dd32$.

\begin{figure}[tbh]
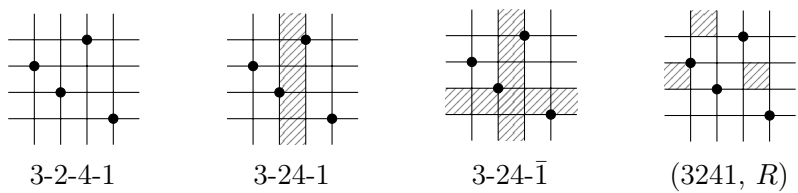

\hfill\hfill
\shortstack{\pattern{}{ 4 }{ 1/3, 2/2, 3/4, 4/1 }{}\\[1ex]3-2-4-1}
\hfill
\shortstack{\pattern{}{ 4 }{ 1/3, 2/2, 3/4, 4/1 }{ 2/0, 2/1, 2/2, 2/3, 2/4 }\\[1ex]3-24-1}
\hfill
\shortstack{\pattern{}{ 4 }{ 1/3, 2/2,  3/4, 4/1 }
            { 2/0, 2/1, 2/2, 2/3, 2/4, 0/1, 1/1, 2/1, 3/1, 4/1 }\\[1ex]3-24-$\bar1$}
\hfill
\shortstack{\pattern{}{ 4 }{ 1/3, 2/2, 3/4, 4/1 }{ 0/2, 1/4, 3/2
  }\\[1ex](3241, $R$)}
\hfill\hfill\hfill
\caption{The patterns 3-2-4-1 (classical), 3-24-1 (vincular),
  3-24-$\bar1$ (bivincular) and the mesh pattern $(3241,R)$, where
  $R=\{(0,2),(1,4),(3,2)\}$.\label{fig-patts}}
\end{figure}

These three kinds of patterns are illustrated in
Figure~\ref{fig-patts}, where the shaded column between the second and
third black dots of the diagram for $3\dd24\dd1$ indicates that in an
occurrence of $3\dd24\dd1$ in a permutation $\pi$ no letter of $\pi$
is allowed to lie between the letters corresponding to the 2 and the 4
or, equivalently, that those letters have to be adjacent in $\pi$.
Patterns thus represented by diagrams with entire columns shaded are
called vincular patterns, but were called generalised patterns when
they were introduced in \cite{babstein}.  Similarly, the shaded row in
the diagram for $p=3\dd24\dd\bar1$ indicates that in an occurrence of
$p$ in a permutation $\pi$, there must be no letters in $\pi$ whose
values lie between those corresponding to the 1 and the 2, that is,
that the letters corresponding to the 1 and 2 must be adjacent in
value, because of the bar over the~1 in $p$.  Patterns with some rows
and columns shaded are called bivincular, and were introduced in
\cite{bcdk}.  

% As an example, the subsequence 142 in 31542 is an occurrence of
% $\bar1\dd32$, since the 4 and 2 are adjacent in position and the 1 and
% 2 are adjacent in value.

Mesh patterns, introduced in \cite{brand-cla-mesh} are now defined by
extending the above prohibitions determined by shaded columns and rows
to a shading of an arbitrary subset of squares in the diagram.  Thus,
in an occurrence, in a permutation $\pi$, of the pattern $(3241,R)$ in
Figure~\ref{fig-patts}, there must, for example, be no letter in $\pi$
that precedes all letters in the occurrence and lies between the
values of those corresponding to the 2 and the~3.  This is required by
the shaded square in the leftmost column.  For example, in the
permutation 415362, 5362 is not an occurrence of $(3241,R)$, since 4
precedes 5 and lies between 5 and 3 in value, whereas the subsequence
4362 is an occurrence of this mesh pattern.

The bivincular patterns were introduced by Bousquet-M\'elou, Claesson,
Dukes and Kitaev in \cite{bcdk}.  The original motivation behind their
definition was to increase the symmetries of the vincular patterns,
whose set is invariant under taking complement and reverse
(corresponding to reflecting their diagrams horizontally and
vertically), but not with respect to taking inverse, which corresponds
to reflecting along the SW-NE diagonal.  The study of the bivincular
pattern $2\dd3\bar1$ was the catalyst of the paper \cite{bcdk}, where
permutations avoiding this pattern were shown to be in bijective
correspondence with two other families of combinatorial objects, the
$(2+2)$-free posets (see Figure~\ref{fig-22p}) and the \emph{ascent
  sequences}.  An ascent sequence is a sequence $a_1a_2\ldots a_n$ of
nonnegative integers where $a_1=0$ and, for all $i$ with $1<i\le n$,
$$ a_i\le\asc(a_1a_2\ldots a_{i-1})+1,
$$ where $\asc(a_1a_2\ldots a_k)$ is the number of \emph{ascents} in
the sequence $a_1a_2\ldots a_k$, that is, the number of places $j\ge1$
such that $a_{j}<a_{j+1}$.  An example of such a sequence is
0101312052, whereas 0012143 is not, because the 4 is greater than
$\asc(00121)+1=3$.  

This connection between $(2+2)$-free posets (also known as
\emph{interval orders}) and ascent sequences led to the determination
of the elegant generating function for these families of objects,
given in \cite[Theorem 13]{bcdk} as
\[
  \sum_{n\ge 0} \ \prod_{i=1}^n \left( 1-(1-t)^i\right).
\]
This generating function, equivalently an exact enumeration of the
$(2+2)$-free posets, had eluded researchers for a long time.  This
became tractable because the previously little studied ascent
sequences are more easily amenable to an effective recursive
decomposition than the posets. The paper \cite{bcdk} has been followed
by a great number of papers on these and bijectively related
combinatorial objects, with no end in sight.  The study of pattern
avoidance by the ascent sequences in their own right was recently
initiated \cite{duncan-es} and has been furthered in \cite{chen-sagan,
  mansour-shattuck-ascent-seqs, yan-ascent-3-nonnesting}, each of
which proves some of the conjectures made in \cite{duncan-es}.

\setlength{\unitlength}{1mm}
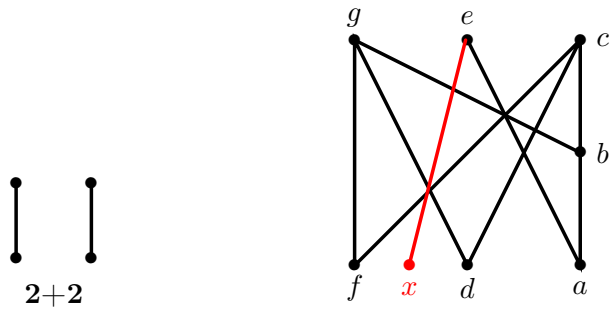
\begin{figure}[tbh]
\newcommand\mbu{\mb{$\bullet$}}
  \def\mb#1{\makebox(0,0){#1}}
  \begin{picture}(100,50)(0,-1)
\allinethickness{.5mm}
\put(30,10){
\path(0,0)(0,10)
\path(10,0)(10,10)
\put(0,0){\mbu}
\put(10,0){\mbu}
\put(0,10){\mbu}
\put(10,10){\mbu}
\put(5,-5){\mb{\bf2+2}}
}
%%%%%%%%%%%%%%%%%%%%%%%%%%%%%%%%%%%%%%%%%%%%%%%%%%%%%%%%%%%%%%%%%%%%%%
\setlength{\unitlength}{3mm}

\put(25,3){
\put(0,0){
\allinethickness{.5mm}
\path(5,10)(10,0)(10,10)(0,0)(0,10)(5,0)(10,10)
\path(0,10)(10,5)
\put(0,0){\mbu}
\put(5,0){\mbu}
\put(10,0){\mbu}
\put(10,5){\mbu}
\put(0,10){\mbu}
\put(5,10){\mbu}
\put(10,10){\mbu}
}
%%%%%%%%%%%%%%%%%%%%%%%%%%%%%%%%%%%%%%%%%%%%%%%%%%%%%%%%%%%%%%%%%%%%%%
\put(0,0){
\put(0,-1){\mb{$f$}}
\put(5,-1){\mb{$d$}}
\put(10,-1){\mb{$a$}}
\put(11,5){\mb{$b$}}
\put(11,10){\mb{$c$}}
\put(5,11){\mb{$e$}}
\put(0,11){\mb{$g$}}
}
%%%%%%%%%%%%%%%%%%%%%%%%%%%%%%%%%%%%%%%%%%%%%%%%%%%%%%%%%%%%%%%%%%%%%%
\put(-.5,0){
\color{red}\path(5,10)(2.5,0)
\put(2.5,-1){\mb{$x$}}
\put(2.5,0){\mbu}
}
}
%%%%%%%%%%%%%%%%%%%%%%%%%%%%%%%%%%%%%%%%%%%%%%%%%%%%%%%%%%%%%%%%%%%%%%
  \end{picture}
  \caption{The poset $\mathbf{2+2}$, consisting of two disjoint
    chains, and a poset containing $\mathbf{2+2}$ in the subposet
    induced by $b,c,e$ and $x$. The poset on the right is
    ($\mathbf{2+2}$)-free if the vertex $x$ is
    removed. \label{fig-22p}}
\end{figure}

We say that two patterns are \emph{Wilf equivalent}, and belong to the
same \emph{Wilf class}, if, for each $n$, the same number of
permutations of length $n$ avoids each.  Of course, the symmetries
mentioned above, corresponding to reflections of the diagram of a
pattern, are (rather trivial) examples of Wilf equivalence, but there
are many more, and they are often hard to prove.  The smallest example
of non-trivial Wilf equivalence is that for the classical patterns of
length 3:\, The patterns $1\dd2\dd3$ and $3\dd2\dd1$ are trivially
Wilf equivalent, and the same is true of the remaining four patterns
of length 3, namely $1\dd3\dd2$,\; $2\dd1\dd3$, $2\dd3\dd1$, and
$3\dd1\dd2$.  All six of these patterns are Wilf equivalent, which is
easy but non-trivial to prove; each is avoided by $C_n$ permutations
of length~$n$, where $C_n$ is the Catalan number
$\frac{1}{n+1}\ch2n,n,$.

By extension, we say that two patterns $p$ and $q$ are \emph{strongly
  Wilf equivalent} if they have the same \emph{distribution} on the
set of permutations of length $n$ for each $n$, that is, if for each
nonnegative integer $k$ the number of permutations of length $n$ with
exactly $k$ occurrences of $p$ is the same as that for $q$. It is easy
to see that the symmetry equivalences mentioned above imply strong
Wilf equivalence.  For example, $p=1\dd3\dd2$ is strongly Wilf
equivalent to $q=2\dd3\dd1$, since the bijection defined by reversing
a permutation turns an occurrence of $p$ into an occurrence of $q$ and
conversely.  On the other hand, $1\dd3\dd2$ and $1\dd2\dd3$ are not
strongly Wilf equivalent, although they are Wilf equivalent.  For
example, the permutation $1234$ has four occurrences of $1\dd2\dd3$,
but there is no permutation of length 4 with four occurrences of
$1\dd3\dd2$,

Just as
finding the distribution of occurrences of a pattern is in general
harder than finding the number of avoiders, so there are still few
results about strong Wilf equivalence.  In fact, the only nontrivial
such results I am aware of are recent results of Kasraoui
\cite{kasraoui}, who gives an infinite family of strong Wilf
equivalences for non-classical vincular patterns, including, for
example, the equivalence of $3\dd421$ and $421\dd3$.

The mesh patterns, which may seem overly general at first sight, turn
out to be the right level of generalisation for expressing a seemingly
deep algebraic relationship (see Section~\ref{sec-algebra}) that
encompasses all the kinds of patterns defined above.  They also allow
for simple expressions of some more cumbersome definitions, such as
some of the so-called barred patterns, and they provide elegant
expressions for some well known statistics on permutations, such as
the number of left-to-right maxima and the number of components in a
permutation $\pi$, that is, the maximum number of terms in a direct
sum decomposition of $\pi$, to be defined in Section~\ref{sec-layered}.

\section{Enumeration and asymptotics of avoiders: The case of
  1324}\label{sec-1324}

For information about the state of the art in enumeration of
permutations avoiding given patterns, refer to \cite{kitaev-book}.
Here I will essentially only treat one unresolved case, that of the
only classical pattern of length 4, up to Wilf equivalence, for which
neither exact enumeration nor asymptotics have been determined.  Of
course, there is an endless list of problems left to solve when it
comes to avoidance, until we find general theorems.  Whether that will
ever happen is likely to remain unknown for quite a while, given the
slow progress so far.  The situation is similar for vincular patterns,
a short survey on which appears in \cite{genpatt-survey}, and even
less has, understandably, been done when it comes to bivincular and
mesh patterns.  Although we are up against a major obstacle in
furthering the knowledge about classical pattern avoidance, there are
probably some reasonably easy, and interesting, results left to be
found for the vincular, bivincular and mesh patterns.

Observe that in the remainder of the paper, unless otherwise noted, I
am talking about classical patterns, which will be written in the
classical way, that is, without any dashes to separate their letters
(as should be done were they being considered as vincular patterns).

In 2004, Marcus and Tardos \cite{marcus-tardos} proved the
Stanley-Wilf conjecture, stating that, for any classical pattern~$p$,
we have $\avpn<C^n$ for some constant $C$ depending only on~$p$, where
$\avpn$ is the number of permutations of length $n$ avoiding $p$.  It
had been shown earlier by Arratia ~\cite{arratia} that this was
equivalent to the existence of the limit
$$
\sw(p)=\lim_{n\rightarrow\infty}{\sqrt[n]{\avpn}},
$$ 
which is called the \emph{Stanley-Wilf limit for $p$}.  It should be
noted that this exponential growth does not apply to vincular patterns
in general.  For example, the avoiders of length $n$ of the pattern
$1\dd23$ are enumerated by the Bell numbers $B_n$ \cite{cla-gpa},
which count set partitions and grow faster than any exponential
function. Namely, when $n$ goes to infinity, the quotient
$B_n/B_{n-1}$ grows like $n/\log n$
\cite[Prop. 2.6]{klazar-set-systems-weight}. In fact, it seems likely
(see \cite[Section 8]{genpatt-survey}) that there are no vincular
patterns of length greater than 3 with exponential growth.

The Stanley-Wilf limit is 4 for all patterns of length 3, which follows from the fact that the number of avoiders of any one of these is the $n$-th Catalan number $C_n$, as mentioned above. This limit is known to be 8 for the pattern $1342$ (see B\'ona's paper \cite{bona-1342}).  For the pattern $1234$ the limit is 9.  This is a special case of a result of Regev \cite{regev} (see also the more recent \cite{regev-asymp-young}), who provided a formula for the asymptotic growth of the number of standard Young tableaux with at most $k$ rows, pairs of which are in bijection, via the Robinson-Schensted correspondence, with permutations avoiding an increasing pattern of length $k+1$.  This limit can also be derived from Gessel's general result \cite{gessel-symmetric} for the number of avoiders of an increasing pattern of any length.

The only Wilf class of patterns of length 4 for which the Stanley-Wilf
limit is unknown is represented by $1324$.  A lower bound of 9.47 was
established by Albert et al. \cite{albert-et-al}, who used an
interesting technique, the \emph{insertion encoding} of a permutation,
introduced in \cite{albert-linton-ruskuc-insertion}.  This encoding
was used to analyse, in an efficient fashion, how a 1324-avoider can
be built up by inserting the letters $1,2,\ldots,n$ in increasing
order.  They then used this to show that a certain subset of
1324-avoiders has insertion encodings that are accepted by a
particular finite automaton, from which they could deduce the lower
bound.

Successively improved upper bounds have been established in several
steps.  The first reasonably small one, was given in
\cite{cla-jel-est}, where 1324-avoiders were shown to inject into the
set of pairs of permutations where one avoids 132 and the other avoids
213, which led to an upper bound of $4\cd4=16$.  Refining the method
used in \cite{cla-jel-est}, B\'ona~\cite{bona-new-1324} was able to
reduce this bound to $7+4\cdot\sqrt3~\approx13.9282$.

In \cite{cla-jel-est} it is conjectured that
$\sw(1324)<e^{\pi\sqrt{2/3}}\approx13.001954$.  This would follow from
another conjecture in \cite{cla-jel-est}, which says that the number of
$1324$-avoiders of length $n$ with a fixed number $k$ of
inversions\footnote{An inversion in a permutation $a_1a_2\ldots a_n$
  is a pair $(i,j)$ such that $i<j$ and $a_i>a_j$.} is increasing as a
function of~$n$.  It is further conjectured that this holds for
avoiders of any classical pattern other than the increasing ones, and
there is some evidence that this in fact applies to all non-increasing
vincular patterns.

Using Markov chain Monte Carlo methods to generate random
$1324$-avoiders, Madras and Liu \cite{madras-liu} estimated that, with
high likelihood, $\sw(1324)$ lies in the interval $[10.71,11.83]$.
Recent computer simulations\footnote{The simulations were done in the
  following way: Using the method of \cite[Sections 2 and
  4]{madras-liu}, I first generated, in each of roughly 155
  independent processes, $10^9$ 1324-avoiders of length 1,000 (the
  initialisation phase in the terminology of \cite{madras-liu}).
  Using each of these 155 generated seed permutations, each in an
  independent process, I then generated a total of
  approx. $2.57\cd10^{12}$ further avoiders, (the data collection
  phase).  Using one in every ten thousand of those avoiders, I found
  the right end of the leftmost occurrence of the pattern 132.  The
  number of that place equals the number of different places where  $n+1=1001$ can be inserted to obtain an avoider of length 1,001 from one of
  length 1,000.  In the limit where $n\rightarrow\infty$ the average
  of the number of this place over all avoiders of a given length~$n$
  would give the S-W limit.  The average I obtained was approximately
  11.01146.  However, the convergence may  be very slow, so, as pointed out to me by Josef Cibulka~\cite{cibulka-personal}, this may be close to the correct value for $n=1,000$ although the actual limit might be significantly greater.}
I have done (unpublished), building on the random generation method of
Madras and Liu, but estimating $\sw(1324)$ in a different way, point
to the actual limit being close to 11.  Given how hard it seems
to determine $\sw(1324)$ makes it understandable that Doron Zeilberger
is claimed to have said \cite{elder:problems-and-co:} that ``Not even
God knows the number of $1324$-avoiders of length 1,000.''  I'm not
sure how good Zeilberger's God is at math, but I believe that some
humans will find this number in the not so distant future.

\section{Are layered patterns the most easily
  avoided?}\label{sec-layered}

A \emph{layered permutation} is a permutation that is a concatenation
of decreasing sequences, each containing smaller letters than in any
of the following sequences.  An example of a layered permutation is
32145768, whose layers are displayed by $321-4-5-76-8$.  A layered
pattern is thus a \emph{direct sum} of decreasing permutations.  The
direct sum $\sigma\oplus\tau$ of two permutations $\sigma$ and $\tau$
is obtained by appending $\tau$ to $\sigma$ after adding the length of
$\sigma$ to each letter of $\tau$.  The \emph{skew sum}
$\sigma\ominus\tau$ of $\sigma$ and~$\tau$ is obtained by prepending
$\sigma$ to $\tau$ after adding the length of $\tau$ to each letter
of~$\sigma$. Thus, for example, we have $3142\oplus231=3142675$ and
$3142\ominus231=6475231$.

Of course, reversing a layered pattern $p$, or taking its complement,
gives a pattern that is Wilf equivalent to $p$, and such a
pattern/permutation might be called \emph{up-layered}, since each
layer is increasing.  Clearly, an up-layered pattern is the skew sum
of increasing permutations. To simplify the discussion in this
section, without changing the traditional definition of layered
permutations, I will abuse notation by letting ``layered'' refer both
to layered and up-layered patterns.

Evidence going back at least to Julian West's thesis
\cite[Section~3.3]{west-thesis} supports the conjecture that among all
patterns of length $k$ the pattern avoided by most permutations of a
sufficiently large length $n$ is a layered pattern. This conjecture
has been around for a long time, and several variations on it, ranging
from asymptotic dominance to the conjecture that $\avsn\le\avtn$ for
all $n$ if $\sigma$ and $\tau$ have the same length and~$\tau$ is
layered but~$\sigma$ isn't.  The latter conjecture here is false, as
pointed out by V\'it Jel\'inek \cite{jelinek-personal}, who noted
that, by \cite[Theorem~4.2]{bona-records}, the non-layered pattern
obtained as the direct sum $p\oplus q$ of the layered pattern $p$ with
layer sizes $(1,2,1,2,\ldots,1,2,1)$ and $q=231$ has a larger
Stanley-Wilf limit than the layered pattern $12\ldots k$, if $k$ is
sufficiently large.  One of the weakest of these conjectures, widely
believed to be true, is made explicit in
~\cite[Conjecture~1]{cla-jel-est} (although it may have been made before)
and says that among all patterns of a given length $k$, the largest
Stanley-Wilf limit is attained by some layered pattern.

The first published instances of such conjectures I am aware of appear
in Burstein's thesis~\cite[Conjectures~9.5]{burstein-thesis}, where
several conjectures are listed, both by Burstein and others, some of
which have been refuted, while others have yet to be proved right or
wrong.  I think these questions are quite interesting, apart from
their intrinsic interest, because any results about them are likely to
be accompanied by quite general results about the enumeration and
asymptotics of pattern avoidance.

Although the evidence is strong in support of the conjecture that the
most easily avoided pattern of any given length is a layered pattern,
there is currently no general conjecture that fits all the known data
about the particular layered patterns with the most avoiders.
However, there are some ideas about what form the most avoided layered
patterns ought to have, and specific conjectures that have not been
shown to be false. 

As is mentioned in Burstein's thesis \cite[see Conjectures 9.10, 9.11
and 9.12]{burstein-thesis} K\'ezdy and Snevily had made some
conjectures about this, and Burstein did too.  In all these cases,
which fall into three classes, depending on the congruence class
modulo~4 of~$k$, the patterns conjectured to have the maximal number
of avoiders have small layers, and are highly symmetric.  For example,
a conjecture of K\'ezdy and Snevily \cite[Conj. 9.10]{burstein-thesis}
that has not been refuted is that for patterns of even length~$k$ the
most avoided pattern is the one with layers
$(1,2,2,\ldots,2,2,1)$. Also, Burstein's Conjecture 9.11 in
\cite{burstein-thesis} that for $k\equiv1\;(\bmod\,4)$ the most
avoided pattern of length $k$ has symmetric layers
$(1,2,\ldots,2,3,2,\ldots,2,1)$ still stands.  Data computed by V\'it
Jel\'inek \cite{jelinek-personal} support these conjectures.
Jel\'inek provided all such data mentioned in this section.  These
computations range over all patterns of lengths up to~10 and
permutations of lengths up to 14, although they are not exhaustive
within these ranges.

For patterns of length $k\equiv3\;(\bmod\,4)$, however, the situation
is not quite that simple.  Burstein's Conjecture 9.12 in
\cite{burstein-thesis}, has the symmetric pattern
(1,2,\ldots,2,1,2,\ldots,2,1) as most avoided, but among permutations
of lengths 11 and~12 the most avoided pattern of length 7 is 1432657,
with layers (1,3,2,1). The next pattern in this respect is 2143657,
with layers (2,2,2,1), and only in third place comes 1324657, with the
symmetric layers (1,2,1,2,1).  The respective numbers of avoiders of
these patterns, for permutations of length 12, are 457657176,
457656206 and 457655768, differing in their last four digits.

But, strange things do occasionally happen in this field, so
Burstein's pretty conjecture for $k\equiv3\;(\bmod\,4)$ cannot be
written off entirely yet as far as asymptotic growth is concerned.
Namely, Stankova and West \cite[Figure~9]{stankova-west} unearthed a
somewhat irregular behaviour in pattern avoidance that gives pause
here.  They computed data that show, among other examples, that
although $\avpn(53241)<\avpn(43251)$ when $7\le n\le12$, the
inequality is switched for $n=13$.  It is still unknown whether such a
switch occurs twice or more for any pair of patterns, but Stankova and
West conjecture \cite[Conjecture~2]{stankova-west} that this does not
happen.  They also make the weaker conjecture that for any pair of
patterns $\sigma$ and~$\tau$ there is an $N$ such that either
$\avsn<\avtn$ for all $n>N$ or else the opposite inequality holds for
all $n>N$.  In other words, that two patterns do not switch places
infinitely often.

It was shown in \cite[Corollary 7]{cla-jel-est} that the Stanley-Wilf
limit of a layered pattern of length $k$ is at most $4k^2$.  In a
recent preprint \cite{bona-best-upper} B\'ona, building on results
from \cite{cla-jel-est} and \cite{bona-best-upper}, proves that this
bound can be reduced to $2.25k^2$.  If the conjecture is true that the
most avoided pattern of length $k$, for any given $k$, is a layered
pattern then this would be an upper bound for the Stanley-Wilf limit
of any pattern of length $k$, but it is not known whether this bound
is the best possible.

% West:  p. 47 

% Alex: The conjectures you're looking for are 9.10, 9.11 and 9.12 on
% page 154. (By the way, there are some other ones there that seem to be
% true, e.g. 9.8 on page 153. Also, the missing piece of data
% |S_5(14325)|=261863 that Julian computed at my request is on page 146.

\section{The pattern poset, its M\"obius function and topology
}\label{sec-poset}

The set of all permutations forms a poset $\clp$ with respect to
classical pattern containment.  That is, a permutation $\sigma$ is
smaller than $\pi$, denoted $\sigma\le\pi$, if $\sigma$ occurs as a
pattern in $\pi$.  This poset is the underlying object of all studies
of pattern avoidance and containment.  In the following two
subsections I treat the M\"obius function of $\clp$, perhaps the most
studied invariant of a poset, and then the topological aspects of the
simplicial complexes naturally associated to intervals in $\clp$.

\subsection{The M\"obius function of the permutation
  poset}\label{sec-mobius}

An \emph{interval} $[x,y]$ in a poset $P$ is the set of all elements
$z\in P$ such that $x\le z\le y$.  An interval is thus a subposet with
unique minimum and maximum elements, namely $x$ and $y$, respectively,
except when $x\not\le y$, in which case $[x,y]$ is empty.  The
\emph{M\"obius function} of an interval is defined recursively as
follows: For all $x$, we set $\mu(x,x)=1$ and
\begin{equation*}
  \mu(x,y)=-\sum_{x\le z<y}{\mu(x,z)}.
\end{equation*}
Thus, the M\"obius function is uniquely defined by setting its sum
over any interval $[x,y]$ to 1 if $x=y$ and to~0 otherwise.
Figure~\ref{fig-mob} gives an example of an interval from $\clp$, and
the values of the M\"obius function on this interval, computed from
bottom to top using the above recursive definition, showing that
$\mu(321,316254)=-1$.

\setlength{\unitlength}{0.35ex}
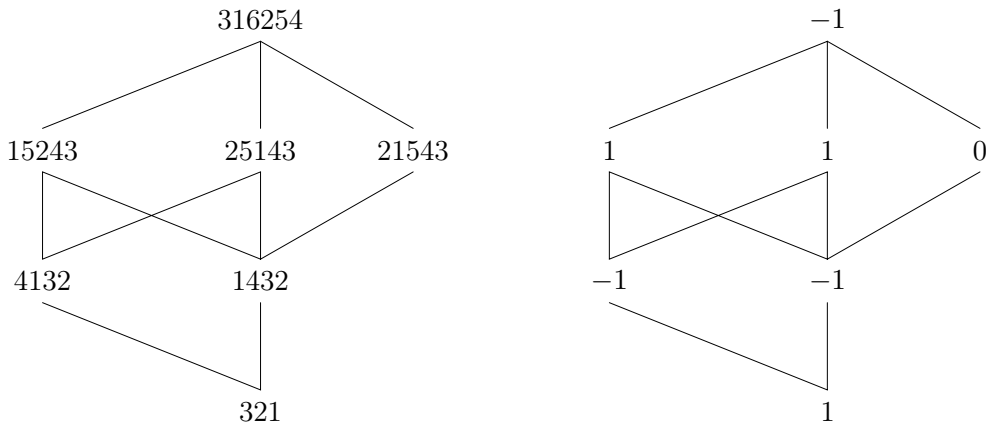
\begin{figure}[htbp]
\begin{picture}(150,110)(-17,-10)
\put(\xb,\yd){\cbox{316254}}
\put(\xa,\yc){\cbox{15243}}
\put(\xb,\yc){\cbox{25143}}
\put(\xc,\yc){\cbox{21543}}

\put(\xa,\yb){\cbox{4132}}
\put(\xb,\yb){\cbox{1432}}

\put(\xb,\ya){\cbox{321}}

\path(\xc,\ycu)(\xb,\ybo)(\xa,\ycu)(\xa,\ybo)(\xb,\ycu)(\xb,\ybo)

\path(\xa,\ybu)(\xb,\yao)(\xb,\ybu)
\path(\xa,\yco)(\xb,\ydu)(\xb,\yco)
\path(\xb,\ydu)(\xc,\yco)

%Right hand picture
\put(130,0){
\put(\xb,\yd){\cbox{$-1$}}
\put(\xa,\yc){\cbox{1}}
\put(\xb,\yc){\cbox{1}}
\put(\xc,\yc){\cbox{0}}

\put(\xa,\yb){\cbox{$-1$}}
\put(\xb,\yb){\cbox{$-1$}}

\put(\xb,\ya){\cbox{1}}

\path(\xc,\ycu)(\xb,\ybo)(\xa,\ycu)(\xa,\ybo)(\xb,\ycu)(\xb,\ybo)

\path(\xa,\ybu)(\xb,\yao)(\xb,\ybu)
\path(\xa,\yco)(\xb,\ydu)(\xb,\yco)
\path(\xb,\ydu)(\xc,\yco)

}
\end{picture}
\caption{The interval $\cli=[321,316254]$ in $\clp$ and the values of
  $\mu(321,\tau)$ on its elements $\tau$. \label{fig-mob}}
\end{figure}

The Möbius function of intervals in the permutation poset exhibits a
great variety in values, even for permutations of small length, and is
seemingly very hard to determine in the general case.  The first such
results were found by Sagan and Vatter \cite{sagan-vatter}, who found
a formula for intervals of layered permutations.

In \cite{ste-tenner} the problem was solved for some special cases,
and in \cite{mob-sep} an effective (polynomial time) formula was given
for all \emph{separable} permutations.  These are the permutations
that can be built from singletons by combinations of direct sums and
skew sums.  The permutation 31254 is separable, since
$31254=(1\om(1\op1))\op(1\om1)$, whereas 3142 is not.  In fact, the
separable permutations are precisely those that avoid both 3142 and
2413 (the two \emph{simple} permutations of length 4; see definition
later in this section).  This fact is usually attributed to folklore,
and it is straightforward to prove.
% , but Vatter gives a proof in lecture notes on the web
% \cite{vatter:lecture-notes}.

The formula in \cite{mob-sep} for the M\"obius function of an interval
of separable permutations $[\sigma,\tau]$ is based on the
representation of separable permutations by rooted trees, and then
certain embeddings of (the tree for) $\sigma$ in (the tree for) $\tau$
are counted, with a sign, to obtain $\mu(\sigma,\tau)$.  More
precisely, the \emph{reduced separating tree} of a separable
permutation $\pi$ describes exactly how $\pi$ is composed by sums and
skew sums.  The separating tree is defined recursively by letting the
children of the root be the summands in the (skew) sum of $\pi$
decpomposed into the maximal number of summands (see
\cite[Figure~1]{mob-sep}).  The embeddings of the tree for $\tau$ into
the tree for $\sigma$ that are counted to compute the M\"obius
function are the so called \emph{normal embeddings}, defined by a
rather technical condition that we won't go into here.  Each embedding
corresponds to a unique occurrence of $\sigma$ in $\tau$.  The formula
for the M\"obius function $\mst$ is then given by
\begin{equation}\label{eq-mf-sep}
\mu(\sigma,\tau) = \sum_{f\in N(\sigma,\tau)}{\sgn(f)},
\end{equation}
where $N(\sigma,\tau)$ is the set of normal embeddings of $\sigma$ in
$\tau$ and $\sgn(f)$ is a sign associated to the embedding $f$,
depending on which leaves of the tree for $\tau$ belong to the
embedding of $\sigma$ in $\tau$.  Again, this is determined by a
rather technical condition that we omit, but we refer the reader to
\cite[Section~4]{mob-sep}.

Formula (\ref{eq-mf-sep}) implies, among many other things, that the
absolute value of $\mu(\sigma,\tau)$ cannot exceed the number of
occurrences of $\sigma$ in $\tau$ (which is far from true in the
general case), because each of the embeddings in question corresponds
to a unique occurrence of $\sigma$ in $\tau$.

This formula also allows for many easy computations of values of the
M\"obius function, such as this example: If
$\pi_i=1,3,5,\dots,2i-1,2i,\dots,4,2$, then
$$
\mu(\pi_k, \pi_n)=\ch n+k-1,n-k,.
$$
For example, $\mu(\pi_2,\pi_4)=\mu(1342,13578642)=\ch4+2-1,4-2,=10$.

It is conjectured in \cite[Conjecture 30]{mob-sep} that the maximum
value of $\mu(\sigma,\pi)$ for any separable $\pi$ of length $n\ge3$
is obtained by a permutation of this form, for $k$ that is roughly
$n/2$ (but whose exact formula depends on the parity of the length of
$\pi$).

In \cite{mob-sep} a recursive formula for the M\"obius function was
also given in the case of \emph{decomposable} permutations, those that
can be written non-trivially as sums or skew sums
($241365=2413\op(1\om1)$ is decomposable, whereas 2413 is not).  This
reduces the problem to the indecomposable permutations, for which it
seems unlikely there will be a general formula anytime soon.  However,
solutions for large classes of these are reasonable to hope for.
Moreover, indications are that the results for separable permutations,
which are characterised by avoiding $2413$ and $3142$, can be extended
to other more complicated permutation classes.  That would allow us to
find the maximum (absolute) value of the M\"obius function on infinite
permutation classes by computing it on a small finite set of
permutations in the class.

An occurrence of a \emph{consecutive} (vincular) pattern $p$ in a
permutation $\pi$ is an occurrence of $p$ in $\pi$ whose letters are
consecutive in $\pi$, such as the occurrence 364 of the pattern 132 in
5136472.  An effective formula for computing the M\"obius function of
the poset of permutations ordered by containment as consecutive
patterns was given in \cite{bern-ferr-ste}.  Sagan and Willenbring
\cite{sagan-willenbring} later provided another proof, and used
discrete Morse theory to determine the homotopy type of intervals of
this poset.  This poset is rather simple; its M\"obius function is
restricted to $0$, $1$ and $-1$, and it was shown in
\cite{sagan-willenbring} that its intervals are either contractible or
else homotopy equivalent to a single sphere (see the next subsection,
on topology of the permutation poset).

A recent paper by McNamara and Sagan \cite{mcnamara-sagan} established
similar and further results in the more wide ranging case of the poset
of generalised subword order.  In particular, they determined the
M\"obius function for all intervals of this poset, and exhibited a
family of intervals that are homotopic to wedges of spheres.  Given
the similarity in definition of this poset to the poset $\clp$ of
permutations with the classical pattern containment order, it seems
reasonable to hope that the methods developed in \cite{mcnamara-sagan}
can be adapted to obtain results for large classes of intervals in
$\clp$.

Let $\one$ be the permutation of length 1. It is shown in
\cite[Corollary 24]{mob-sep} that for any separable permutation $\pi$
the only possible values of $\mu(\one,\pi)$ are 0, 1 and $-1$.
Although nobody has bothered proving this yet, it appears certain that
for arbitrary permutations $\pi$ the absolute value $\mu(\one,\pi)$ is
unbounded.  For example, it seems a safe guess (based on computed
data) that
$$
\mu(\one, 2468\ldots(2n)135\ldots(2n-1)) = -\ch n+1, 2,.
$$
As an example, $\mu(\one,24681357)=-\ch4+1,2,=-10$.  

The maximum absolute value of $\mu(\one,\pi)$ is known for all $\pi$
of length at most $11$, as mentioned in \cite[Section~5]{mob-sep}.
Recent computations I have done, although not exhaustive, suggest that
for $n=12$ the maximum is attained only by the permutation
$\pi=4\;7\;2\;10\;5\;1\;12\;8\;3\;11\;6\;9$, and its symmetric
equivalents, with $\mu(\one,\pi)=-261$.  That would match the results
for $n<12$, mentioned in \cite{mob-sep}, namely that the maximum is in
each case attained by only one permutation (up to trivial symmetries),
and that that permutation is \emph{simple}\footnote{Except when $n=3$,
  for which there are no simple permutations}.  A permutation
$\pi=a_1a_2\ldots a_n$ is simple if it has no segment
$a_ia_{i+1}\ldots a_{i+k}$, where $0<k<n-1$, that consists of a
segment of values, that is, such that $\{a_i,a_{i+1},\ldots,
a_{i+k}\}=\{\ell,\ell+1,\ldots,\ell+k\}$ for some $\ell$.  For
example, 315264 is simple, but 461325 is not, since 132 is a segment
of consecutive values that is nontrivial, that is, neither a singleton
nor the entire permutation.  A nice survey on simple permutations,
which are important in the study of permutation classes, is found in
\cite{brignall:a-survey-of-sim:}.

It is worth noting here that separable permutations and simple
permutations are in some imprecise sense each others' opposites; the
separable ones decompose very nicely, while being simple is an
obstruction to such decomposition.  Thus, it would not be surprising
if the maximum of $|\mu(\one,\pi)|$ over all $\pi$ of length $n$ turns
out to be attained by a simple permutation.  It seems less certain
that there is, for all $n$, a unique permutation, up to trivial
symmetries, that attains the maximum value for each $n$.  The sequence
of values of $\mu(\one,\pi)$ for which $|\mu(\one,\pi)|$ is maximised
(for each length $n$, starting at $n=1$) begins with
$$
1, -1, 1, -3, 6, -11, 15, -27, -50, -58, 143, -261, \ldots
$$
(the last entry still conjectural, as mentioned above) but no
nontrivial upper bound on its $n$-th term is known.

Although many families of intervals $\ist$ with $\mst=0$ are described
in \cite{mob-sep,ste-tenner}, it is an open problem to characterise
such intervals completely.  It might also be interesting to
characterise those intervals for which $|\mu(\sigma,\pi)|$ equals the
number of occurrences of $\sigma$ in $\pi$.

Another question raised in \cite{mob-sep} is whether it is possible to
find a bound on $|\mst|$ that depends only on the number of
occurrences of $\sigma$ in $\tau$.  As mentioned above, it has been
shown that for separable $\sigma$ and $\tau$, the value $|\mst|$
cannot exceed the number of occurrences of $\sigma$ in $\tau$, but
nothing similar is known for the general case.

In conclusion, although a general formula for the M\"obius function of
an interval will likely remain untractable for a while it seems
reasonable to expect much further progress.  In particular, the fact
that there are families of intervals whose M\"obius function has a
``nice'' formula, such as binomial coefficients, gives hope that these
intervals have a structure that can be understood and exploited, and
used to elicit the topology of these intervals, which is treated in
the next section.

\subsection{Topology of the permutation poset}\label{sec-topo}

Another important aspect of the poset of permutations, as for any
combinatorially defined poset, is the topology of (the \emph{order
  complexes} of) its intervals.  (For terminology not defined in this
section, see \cite{cca}.)  The order complex of a poset $P$, denoted
$\Delta(P)$, is the abstract simplicial complex consisting of the
\emph{chains} of $P$.  A chain in a poset $P$ is a set of elements in
$P$ that are pairwise comparable, and thus totally ordered.  Any
subset of a chain forms a chain, so the set of chains is closed in
this respect, which is a defining property of simplicial complexes.
To study the topology of an interval $\cli$ we remove the maximum and
minimum element of $\cli$, and take the order complex of the remaining
``interior'' of $\cli$, which we denote by $\bar\cli$.

There are some indications that large classes of intervals in the
permutation poset $\clp$ may have a ``nice'' topology, meaning that
the topology can be simply described, and thus that the homology, and
the homotopy type, of these intervals can be well understood.  Gaining
such understanding may well lead to significant progress in answering
other questions, for example about the M\"obius function.  This is
because, in some sense, the topological properties of a complicated
interval give a much clearer picture of the important overall
structure, sweeping aside irrelevant details that obscure the view.
One well known connection to the M\"obius function is that the
M\"obius function of a poset is equal to the \emph{reduced Euler
  characteristic} of its order complex, which is a topological
invariant (see \cite[Proposition 3.8.6]{ECI}).

Figure~\ref{fig-cplx} shows the interval $\cli=[321,316254]$, and the
order complex of $\bar\cli$.  This complex is homotopy equivalent to a
sphere, since contracting the edge between 1432 and 21543 leaves a
1-dimensional sphere (topologically speaking) consisting of the edges
of the rectangle in the figure.  Note that $\mu(321,21543)=0$ and so
removing 21543 from $\cli$ does not affect the M\"obius function of
$\cli$.  The reduced Euler characteristic of a 1-dimensional sphere is
$-1$, which, of course, equals $\mu(321,316254)$.  An obvious question
(see below) is whether it is common for order complexes arising in
this way from intervals of $\clp$ to have similarly nice properties,
such as being homotopy equivalent to wedges of spheres of the same
dimension.  

The elements (which are sets) in a simplicial complex are called
\emph{faces}.  A \emph{facet} of a simplicial complex is a face that
is maximal with respect to containment.  A simplicial complex is
\emph{pure} if all its facets have the same dimension.  A property
that implies a pure simplicial complex is homotopy equivalent to a
wedge of spheres is being \emph{shellable}.  Informally, this means
that the complex can be built up, one facet at a time, such that each
facet that is added, apart from the first one, intersects the union of
the previous ones in a pure subcomplex of the maximum possible
dimension (which is one less than the dimension of the facets).  The
complex in Figure~\ref{fig-cplx} is shellable.  One shelling order of
its facets (edges) is $a,b,c,d,e$, whereas beginning with $a,c,\ldots$
can not give a shelling since the edge $c$ does not intersect $a$.  As
mentioned before, a shellable complex is necessarily homotopy
equivalent to a wedge of spheres, and showing shellability is probably
the most common tool used to determine the topology of combinatorially
defined complexes.

\setlength{\unitlength}{0.35ex}
\begin{figure}[btp]
\begin{picture}(150,110)(-17,-10)
\put(\xb,\yd){\cbox{316254}}
\put(\xa,\yc){\cbox{15243}}
\put(\xb,\yc){\cbox{25143}}
\put(\xc,\yc){\cbox{21543}}

\put(\xa,\yb){\cbox{4132}}
\put(\xb,\yb){\cbox{1432}}

\put(\xb,\ya){\cbox{321}}

\path(\xc,\ycu)(\xb,\ybo)(\xa,\ycu)(\xa,\ybo)(\xb,\ycu)(\xb,\ybo)

\path(\xa,\ybu)(\xb,\yao)(\xb,\ybu)
\path(\xa,\yco)(\xb,\ydu)(\xb,\yco)
\path(\xb,\ydu)(\xc,\yco)

\put(130,0){
\def\xa{0}
\def\xb{30}
\def\xc{65}
\def\xaa{15}
\def\yaa{10}
\def\yba{15}

\put(\xa,\yc){\cbox{15243}}
\put(\xb,\yc){\cbox{1432}}
\put(\xc,\yc){\cbox{21543}}

\put(\xa,\yb){\cbox{4132}}
\put(\xb,\yb){\cbox{25143}}

\put(\xaa,\yaa){\cbox{ $\Delta(\bar\cli)$}}

\path(\xc,\ycu)(\xa,\ycu)(\xa,\ybo)(\xb,\ybo)(\xb,\ycu)

\put(\xc,\ycu){\cbox{\p}}

\put(\xb,\ybo){\cbox{\p}}
\put(\xa,\ybo){\cbox{\p}}
\put(\xa,\ycu){\cbox{\p}}
\put(\xb,\ycu){\cbox{\p}}

\put(\xaa,38){\cbox{$a$}}
\put(4,44.5){\cbox{$b$}}
\put(\xaa,51){\cbox{$c$}}
\put(27,44.5){\cbox{$d$}}
\put(45,51){\cbox{$e$}}

}
\end{picture}
\caption{The interval $\cli=[321,316254]$ in $\clp$ and the order
  complex of $\bar\cli$.  Note that this is a deceptively simple
  example, since intervals of a higher rank have order complexes of
  higher dimension, which are not so easy to depict.\label{fig-cplx}}
\end{figure}
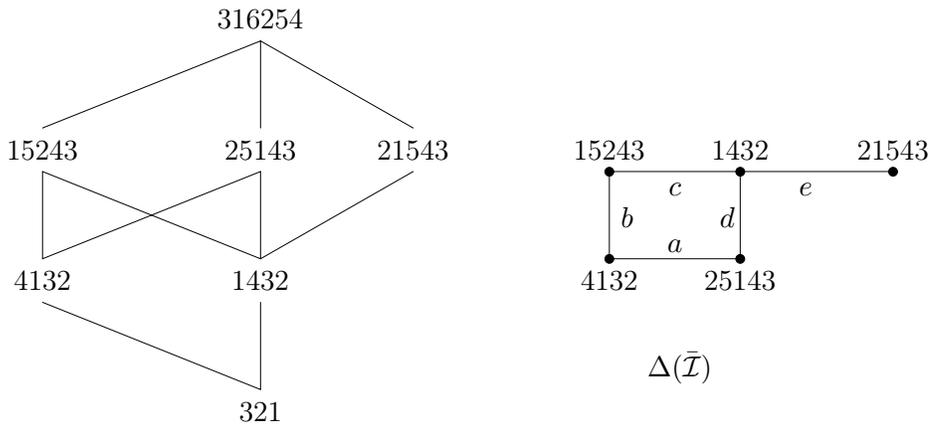

As mentioned in the previous section, Sagan and Willenbring
\cite{sagan-willenbring} used discrete Morse theory to determine the
homotopy type of intervals of the poset of consecutive patterns, whose
intervals are either contractible or else homotopy equivalent to a
single sphere, and it is reasonable to expect that results using
similar techniques will show large classes of intervals in $\clp$ to
be homotopy equivalent to wedges of spheres.  The M\"obius function
of such intervals is, up to a sign, just the number of spheres in the
wedge.

\def\dsp{\ensuremath{\Delta(\sigma,\pi)}}

In \cite[Question 31]{mob-sep} the following questions were raised,
where, given an interval $\cli=[\sigma,\pi]\in\clp$, we let $\dsp$ be the
order complex of $\bar\cli$:

\begin{enumerate}
\item For which $\sigma$ and $\pi$ does $\dsp$ have the homotopy type
  of a wedge of spheres?
\item\label{q-2} Let $\Gamma$ be the subcomplex of $\dsp$ induced by
  those elements $\tau$ of $[\sigma,\pi]$ for which
  $\mu(\sigma,\tau)\ne0$.  Is $\Gamma$ a pure complex, that is, do all
  its maximal simplices (with respect to inclusion) have the same
  dimension?
\item \label{q-3} If $\sigma$ occurs precisely once in $\pi$, and
  $\mu(\sigma,\pi)=\pm1$, is $\dsp$ homotopy equivalent to a sphere?
\item\label{q-4} For which $\sigma$ and $\pi$ is $\dsp$ shellable?
\end{enumerate}

An example where $\dsp$ is not shellable is given by $\sigma=231$ and
$\pi=231564$, since $\dsp$ in this case consists of two disjoint
components, each of which is contractible.  If, however, we remove
from $[231,231564]$ all those elements $\tau$ for which
$\mu(231,\tau)=0$, we get a shellable complex, namely a boolean
algebra of rank 2.  For questions \ref{q-2} and~\ref{q-3} above we
don't know any counterexamples.  In fact, we know no counterexamples
to question~\ref{q-3} even without the condition of just one
occurrence.  However, since we have so far only examined intervals of
small rank, our evidence is weak.

\section{Other algebraic aspects}\label{sec-algebra}

In addition to the M\"obius function (and the underlying incidence
algebra of $\clp$) there are some algebraic aspects of permutation
patterns that have been little studied, but which might harbour some
interesting things.  I mention two here: The ring of functions of
vincular patterns, and Br\"and\'en and Claesson's reciprocity theorem
for mesh patterns.

Vincular patterns can be regarded as functions from the set of all
permutations to the ring of integers, counting occurrences of
themselves in a permutation.  For example, $2\dd31(416253)=2$,
corresponding to the 462 and 453, which are all the occurrences of
$2\dd31$ in 416253.  Linear combinations of vincular patterns were
used in \cite{babstein} to classify the Mahonian permutation
statistics, which are those that are equidistributed with the number
of inversions.  Other such combinations played a crucial role in
\cite{stein-williams}, where they were used to record the distribution
of various statistics on the filled Young tableaux treated there.

No further work seems to have been done along these lines, although it
is almost certain that there are many equivalences to be found of the
kind discussed in \cite{babstein}.  A promising indication in this
context is that computer experiments suggest that the distribution of
two linear combinations of patterns is the same for all $n$ provided
that it is the same for all $n$ smaller than some (small) constant
depending only on the length of the patterns involved.  As an example,
exhaustive computer search shows that if two linear combinations of
three vincular patterns of length 3 diverge for $n<10$, then they
diverge already for $n=6$.  For vincular patterns of length 4 all such
combinations with different distributions for $n=9$ differ already for
$n=8$.  Of course, it is possible that divergence will occur again for
greater $n$, but this seems unlikely.  A general theorem to this
effect, guaranteeing that checking such equidistributions for small
$n$ is sufficient to establish equidistribution for all $n$, would be
a major breakthrough, provided the values that need to be checked are
small enough, since that would give automatic proofs of various
theorems.  More importantly, proving such a theorem would undoubtedly
require a general understanding we lack today, and thus likely lead to
significant other progress.
  
Seen as functions, as described above, the set of vincular patterns
constitutes a ring of functions.  (It is a tedious but straightforward
exercise to verify that the product of two vincular patterns can be
expressed as a sum of vincular patterns.)  One relation is known in
this ring, namely the \emph{upgrading} mentioned in
\cite[Equation~(2)]{babstein}, an example of which is 
\begin{eqnarray*}
  (21\dd3) = (21\dd43) + (21\dd34) + (31\dd24) + (32\dd14) + (213).
\end{eqnarray*}    
Since this ring contains all linear combinations of vincular patterns,
it might be worthwhile to study its algebraic structure further.  In
particular, it would be interesting to know if there are other
relations in this ring.

The set of mesh patterns also forms a ring of functions that should be
further investigated for its properties.  It is in this ring that the
striking Reciprocity Theorem of Br\"and\'en and Claesson lives
\cite{brand-cla-mesh}.  The Reciprocity Theorem expresses any mesh
pattern (including the classical patterns) as a (possibly infinite)
linear combination of classical patterns whose coefficients are
obtained from values of the \emph{dual} pattern on permutations.  To
express that theorem, let $p=(\pi,R)$ be a mesh pattern, and let $R^c$
be the complement of $R$, that is, $R^c=[0,n]^2\setminus R$, where $n$
is the length of $\pi$.  We then define the dual pattern of $p$ as
$p^*=(\pi,R^c)$, and define $\lambda(\sigma)$ by
$\lambda(\sigma)=(-1)^{n-k}p^*(\sigma)$, where $k$ is the length of
$\sigma$.  The Reciprocity Theorem is then the following identity,
where the sum is over all classical patterns $\sigma$:
\begin{eqnarray*}
p=\sum_{\sigma\in\cls}\lambda(\sigma)\sigma.
\end{eqnarray*}

It seems likely that much can be gained from this theorem, due to its
universal nature.

\section{Growth rates of permutation classes}\label{sec-growth}
 
A \emph{permutation class} is a set of permutations that is closed
with respect to containment.  That is, if $\pi\in\clc$ for a class
$\clc$, and $\sigma$ occurs as a pattern in $\pi$, then
$\sigma\in\clc$.  The set of permutations avoiding any classical
pattern, or set of such patterns, is easily seen to be a class (which
is \emph{not} true for vincular, bivincular or mesh patterns) and
every permutation class is characterised by a unique antichain of
permutations that are avoided by all elements of $\clc$.  That
antichain is called the \emph{basis} of $\clc$.  Note that there are
infinite antichains of permutations (see
Brignall~\cite{brignall:grid-classes-an:} for the most general
construction to-date), so bases can be infinite.

Studies of the poset of permutations have in recent years yielded many
results about the diverse collection of permutation classes,
parallelling work being done on other types of object (surveyed in
Bollob\'as~\cite{bollobas:hereditary-and-}).  One of the most active
and successful avenues of investigation has been into the \emph{growth
  rates} of permutation classes.  Given a class $\clc$, where $\cln$
is a set of permutations in $\clc$ of length $n$, the growth rate of
$\clc$ is defined as
$$
\textrm{gr}(\clc)=\limsup_{n\rightarrow\infty}\sqrt[n]{|\cln|}.
$$
To connect this terminology with that of Section~\ref{sec-1324}, note
that the Stanley-Wilf limit of the (classical) pattern $p$ is the
growth rate of the class of $p$-avoiding permutations.  Thus the
Stanley-Wilf Conjecture in this context states that all proper
permutation classes have finite growth rates.  One of the most natural
open questions is whether the limit superior above in the definition
can be replaced by a limit; this is known to be possible in the case
of singleton-based classes by Arratia~\cite{arratia}.

The line of research on growth rates attempts to determine both which
growth rates are possible and where notable phase transitions take
place in this spectrum.  The first answers were provided by Kaiser and
Klazar~\cite{kaiser:on-growth-rates:}, who characterised the growth
rates up to $2$.  At the smallest end of the scale, it is clear that
$0$ and $1$ are growth rates of permutation classes, and that no
classes have growth rates between these two numbers.  Kaiser and
Klazar showed that the next growth rate is the golden ratio and
established the stronger result that if $|\cln|<F_n$ (the $n$th
Fibonacci number) for \emph{any} $n$, then $|\cln|$ is eventually
polynomial (the structural properties of such classes were later
explored by Huczynska and Vatter~\cite{huczynska:grid-classes-an:}).
Between the golden ratio and $2$, Kaiser and Klazar showed that all
growth rates are roots of $x^k-x^{k-1}-\cdots-x-1$ for some $k$; note
that this makes $2$ the least accumulation point of growth rates.

Vatter~\cite{vatter:small-permutati:} extended the characterisation of
growth rates up to $\kappa\approx 2.21$, the unique positive root of
$x^3-2x^2-1$.  Moreover, $\kappa$ represents a sharp phase transition:
There are only countably many permutation classes of growth rate less
than $\kappa$, but because infinite antichains of permutations begin
to appear at this growth rate, there are uncountably many permutation
classes of growth rate $\kappa$.  Viewed on the number-line of growth
rates, $\kappa$ also lies in an interesting place, as it is the least
accumulation point of accumulation points of growth rates.  Recent
work by Albert, Ru\v{s}kuc, and Vatter~\cite{albert:inflations-of-g:}
has established another threshold at $\kappa$: Every permutation class
of growth rate less than $\kappa$ has a rational generating function,
while there are (by an elementary counting argument using the
existence of infinite antichains) permutation classes of growth rate
$\kappa$ whose generating functions are not even holonomic.  (A
function on the natural numbers is holonomic if it satisfies a linear
homogeneous recurrence relation with polynomial coefficients.)

Working in the closely related context of ordered graphs, Balogh,
Bollob\'as, and Morris~\cite{balogh:hereditary-prop:ordgraphs}
characterised the growth rates up to $2$ and made two conjectures
which would have implied that growth rates of permutation classes are
always algebraic integers and that the set of growth rates contains no
accumulation points from above.  These conjectures were both disproved
by Albert and Linton~\cite{albert:growing-at-a-pe:}, who constructed
an uncountable set of growth rates.  It remains open if, as suggested
by Klazar~\cite{klazar:overview-of-som}, the conjectures of Balogh,
Bollob\'as, and Morris hold when restricted to finitely based
permutation classes.

Building on the work of Albert and Linton,
Vatter~\cite{vatter:permutation-cla} showed that every real number
greater than or equal to $\lambda\approx 2.48$, the unique positive
root of $x^5-2x^4-2x^2-2x-1$, is the growth rate of a permutation
class, and conjectured that $\lambda$ is best possible.

Thus a striking problem remains: To characterise the growth rates
between $\kappa$ and $\lambda$.  While it may well be impossible to
describe the set of growth rates once it becomes uncountable (and
before it consists of all real numbers), one could perhaps hope to
describe it up to this point.  The work of
Vatter~\cite{vatter:permutation-cla} implies that this happens at or
before $\xi\approx 2.32$, the unique positive root of
$x^5-2x^4-x^2-x-1$.  Thus, just as $\kappa$ represents the transition
from countably many to uncountably many permutation classes, $\xi$ may
represent the transition from countably many to uncountably many
growth rates.

\section{Acknowledgements}

I am deeply grateful to Vince Vatter, who provided invaluable help, in particular with the section on growth rates.  I thank Alex Burstein, Mikl\'os B\'ona, V\'it Jel\'inek and Josef Cibulka for helpful information about the subject of Section~\ref{sec-layered}. Jel\'inek also generously provided all the computations mentioned in that section.  I am also indebted to a meticulous and insightful referee who made many good suggestions, pointed out several errors and spurred me to make the paper somewhat more extensive, all of which has improved it significantly from its original version.

\bibliographystyle{abbrv}
\bibliography{allrefs}

\end{document}